\documentclass[11pt,a4paper,oneside]{article}
\usepackage[english]{babel}
\usepackage[latin1]{inputenc}
\usepackage{enumerate}
\usepackage{amssymb}
\usepackage{amscd}
\def\haken#1{\underline{#1}{\raise -0.3ex\hbox{\vphantom{$#1$}\vrule height 0.7ex}}}

\begin{document}

\begin{center}
{\huge \bf About Cartan-Subalgebras in
Lie-Algebras associated to Associative Algebras}
\end{center}

\begin{center}
{\Large Dr. Sven Wirsing}\\ {\Large Gertbergstr. 23}\\ {\Large 69437
Neckargerach}\\ {\Large Sven.Wirsing@Brandt-Partner.de}
\end{center}
\begin{center}
{\Large For Lena}
\end{center}

\section*{Abstract}
Every associative $K$-algebra $A$ is with respect to the multiplication
$a\circ b := ab-ba$ for all $a,b\in A$ a
Lie-Algebra $A^{\circ}$, also known as the Lie-Algebra associated with $A$. In
\cite{ss} S. Siciliano studies Cartan-Subalgebras of
$A^{\circ}$. These are nilpotent subalgebras $C$ of $A^{\circ}$ which
coincide with the normalizer $N_{A^{\circ}}(C)$ of $C$ in $A$. \\ \\For
every finite-dimensional associative unitary $K$-algebra $A$
Siciliano proofs that the Cartan-Subalgebras of
$A^{\circ}$ are exactly the centralizers of the maximal tori of $A$.
Especially, Cartan-Subalgebras of $A^{\circ}$ are subalgebras of the associative algebra $A$.
A torus is a commutative unitary subalgebra of $A$ for which every element
is separable over $K$. An element is separable if its minimal polynomial
is a product of pairwise irreducible separable polynomials of $K[t]$.
For the algebra classes "finite-dimensional central division algebras" and
"finite-dimensional soluble algebras" S.Siciliano describes the Cartan-Subalgebras
in terms of 'maximal separable subfields' and 'radical complements'.
\\\\
At first we give alternative proves for Sicilianos results concerning these two algebra classes. As an example
we compute the Cartan-Subalgebras for soluble group algebras and especially for soluble group algebras related to
dihedral groups.\\
After this we describe the Cartan-Subalgebras for finite-dimensional associative
division, simple, semi-simple and reduced algebras.\\
As a associative subalgebra we analyze the associative structure of
Cartan-Subalgebras. This investigation is closely connected to the
question for which finite-dimensional associative unitary algebras its group
of units is nilpotent.\\
We close this article by giving a strategy for computing Cartan-Subalgebras for associative
algebras with separable radical complements. An easy consequence is again the description of the Cartan-Subalgebras
in the case of soluble algebras. We demonstrate this strategy for group algebras related to dihedral groups.

\section*{Cartan-Subalgebras are associative subalgebras}
S.Siciliano especially proves (see \cite{ss}, theorem 1) that
for every finite-dimensional associative unitary $K$-algebra every
Cartan-Subalgebra of its associated Lie-Algebra is an associative subalgebra.
We generalize this statement to arbitrary associative algebras.

\subsection*{Definition 1}

\begin{enumerate}[(i)]

 \item For all $n\in \mathbb{N}$ we set $\haken{n}:=\mathbb{N}_{\le n}$.

 \item For all $n\in \mathbb{N}_{\ge 2}$ we define\\
       $T_n:= \{(\alpha,\beta) \mid \exists r\in \haken{n-1},
       a_1,...,a_n \in \haken{n}: \alpha = (a_1,...,a_r),$\\
       $\beta = (a_{r+1},...,a_n), a_1<...< a_r, a_{r+1}<...< a_n,
       \haken{n} = \{a_1,...,a_n\}\}$.

 \item For every Lie-Algebra $L$ and $l\in L$ let $ad(l)$ be the
       adjoint representation of $l$ with respect to $L$.

 \item For every nilpotent Lie-algebra $L$ let $cl(L)$ be the
       nilpotency class of $L$.

\end{enumerate}

\subsection*{Remark 1} Let $A$ be an associative $K$-algebra and $a,b\in A$.

\begin{enumerate}[(i)]

 \item For all $h_1\in A$ the derivation rule\\
       $(ab)\,ad(h_1) = (a\,ad(h_1))b + a(b\,ad(h_1))$ is valid.

 \item For all $n\in \mathbb{N}_{\ge 2}, h_1,...,h_n \in A$
       we get by (i) and induction
       {\begin{eqnarray*}
       (ab)ad(h_1)...ad(h_n) & = &(a\,ad(h_1)...ad(h_n))b+a(b\,ad(h_1)...ad(h_n))\\
        & + & \sum\limits_{((\alpha_1,...,\alpha_r),(\alpha_{r+1},...,\alpha_n))\in T_n}
       (a\,ad(h_{\alpha_1})...ad(h_{\alpha_r}))
       (b\,ad(h_{\alpha_{r+1}})...ad(h_{\alpha_n})).\diamond
       \end{eqnarray*}}
\end{enumerate}

\subsection*{Proposition 1}
 Let $A$ be an associative (unitary) $K$-algebra and $C$ be a Cartan-Subalgebra of
 $A^{\circ}$. Then $C$ is an (unitary) associative subalgebra of $A$.\\
 \\
 {\bf\underline{Proof:}}As a subalgebra of $A^{\circ}$ the set $C$ is a $K$-subspace
 of $A$. Let $a,b\in C$. Then we have to prove $ab \in C=N_{A^{\circ}}(C)$.
 This is equivalent to $(ab)ad(h_1) \in C = N_{A^{\circ}}(C)$ for all
 $h_1\in C$. By induction, we must show that there exists a $n\in \mathbb{N}$ such that
 for all $h_1,...,h_n \in C$ the statement
 $(ab)ad(h_1)...ad(h_n)\in C$ is valid. Let
 $n:=2\cdot \, cl(C)$. Using part (ii) of remark 1 the statement
 $(ab)ad(h_1)...ad(h_n)= 0 \in C$ is valid (At least one factor of every summand in the sum of part (ii) in remark 1 is equal to zero.).
 Is $A$ unitary then $1_A \in N_{A^{\circ}}(C)=C$. $\diamond$

\section*{Division algebras}
S. Siciliano proves by theorem 1 in \cite{ss} that the separable maximal subfields of a
finite-dimensional associative central $K$-division algebra are exactly
the Cartan-Subalgebras of its associated Lie-algebra.
Here we give an alternative proof of his theorem and describe the Cartan-Subalgebras for
not necessary central division algebras.

\subsection*{Remark 2} Let $A,B$ be associative $K$-algebras.

\begin{enumerate}[(i)]

 \item For all $a_1, a_2 \in A, b_1, b_2 \in B$ the equation\\
       $(a_1\otimes b_1) \circ (a_2\otimes b_2)=(a_1\circ a_2)
       \otimes (b_1b_2) + (a_1a_2)\otimes (b_1\circ b_2)$ is valid.

 \item Let $A^{\circ}$ be nilpotent and $B^{\circ}$ be abelian. By (i) and induction 
 the algebra $(A\otimes B)^{\circ}$ is nilpotent with the same nilpotency class as $A^{\circ}$.$\diamond$

\end{enumerate}

\subsection*{Remark 3} Let $K$ be a field and $n\in
\mathbb{N}$ with $n\ge 2$. Then $gl(n,K)$ is not nilpotent.$\diamond$\\
\\
Let $A$ be an associative $K$-algebra and $T,S$ be subsets of $A$.
Then we call $C_A(T)$ the centralizer of $T$ in $A$, define $C_S(T):=C_A(T)\cap
S$ and call $Z(A):=C_A(A)$ the center of $A$. Additionally
for all $n\in \mathbb{N}$ we call $A^{n\times n}$ the
algebra of all $n\times n$-matrices over $A$. For a
central-simple finite-dimensional associative $K$-Algebra let
$ind(D)(=ind_K(A))$ be the index of $D$.

\subsection*{Proposition 2} Let $D$ be a finite-dimensional
associative non-commutative $K$-division algebra. Then
$D^{\circ}$ is not nilpotent.\\
\\
{\bf\underline{Proof:}} $D$ is as $Z(D)$-algebra a central
division algebra, and $D^{\circ}$ as $K$-algebra nilpotent if and only if
$D^{\circ}$ is nilpotent as $Z(D)$-algebra. Therefor we assume that $D$ is central.\\
Let $T$ be a maximal subfield of $D^{\circ}$, and we assume that $D^{\circ}$ is nilpotent.
Then $(D\otimes T)^{\circ}$ is nilpotent by remark 2.
It is well-known (see e.g. \cite{rp}) that $D\otimes
T$ and $T^{ind(D)\times ind(D)}$ are isomorphic. By using
remark 3 we get $n=1$.$\diamond$\\
\\
Now we can prove the following enhancement of theorem 2 in \cite{ss} and of a
theorem by E. Noether:

\subsection*{Theorem 1}
Let $D$ be a finite-dimensional associative central $K$-division algebra.

\begin{enumerate}[(i)]

 \item The maximal separable subfields of $D$ are exactly
       the separable maximal subfields of $D$.

 \item There exists a separable maximal subfield. (Noether)

 \item The Cartan-Subalgebras of $D^{\circ}$ are exactly
       the separable maximal subfields of $D$. (Siciliano)

\end{enumerate}

{\bf\underline{Proof:}} ad(i): Let $T$ be a maximal separable
subfield of $D$. Then $T$ is a maximal torus of $D$ (Every unitary subalgebra of $D$ 
is a division algebra and tori are commutative.). By Theorem 1 in \cite{ss} the algebra $C_D(T)$ is a
Cartan-Subalgebra of $D^{\circ}$. As $D$
is a finite-dimensional associative $K$-division algebra so is
$C_D(T)$ as well, and by Proposition 2 the algebra $C_D(T)$ is a subfield of $D$.
$C_D(T)$ is maximal Lie-nilpotent so that $C_D(T)$ is a maximal
subfield of $D$. Maximal subfields are self-centralizing (see for instance
\cite{rp}) so that $C_D(C_D(T))=C_D(T)$ is valid. By the double-centralizer-theorem
we get $C_D(C_D(T))=T$. Thus $T=C_D(T)$ is a separable maximal subfield of $D$. Obviously every separable maximal 
subfield is a maximal separable subfield of $D$.\\
\\
ad(ii): $K\cdot 1_D$ is a separable subfield of $D$. Thus there exists a
maximal separable subfield of $D$ which is by (i) a
separable maximal subfield of $D$. \\
\\
ad(iii): By theorem 1 in \cite{ss} the Cartan-Subalgebras
of $D^{\circ}$ are exactly the centralizers of the maximal tori of
$D$. A maximal torus of $D$ is a maximal separable subfield of $D$ (Every unitary subalgebra of $D$
is a $K$-division subalgebra of $D$ and a torus is commutative.).
By (i) we get that the Cartan-Subalgebras of $D^{\circ}$ are exactly
the centralizers of the separable maximal subfields of $D$.
Now the proof is finished because (see for instance \cite{rp})
every maximal subfield of $D$ is self-centralizing. $\diamond$\\
\\
A consequence of theorem 1 is:

\subsection*{Corollary 1}
Let $D$ be a finite-dimensional associative central $K$-division algebra.

\begin{enumerate}[(i)]

 \item The Cartan-Subalgebras of $D^{\circ}$ are Lie-isomorphic and
       $ind(D)$-dimensional.

 \item For a perfect field $K$ all Cartan-Subalgebras of
       $D^{\circ}$ are exactly the maximal subfields of $D$.$\diamond$

\end{enumerate}

We extend theorem 1 to finite-dimensional associative and not
necessary central $K$-division algebras:

\subsection*{Theorem 2}
Let $D$ be a finite-dimensional associative $K$-division algebra.

\begin{enumerate}[(i)]

 \item The Cartan-Subalgebras of $D^{\circ}$ are exactly
       the maximal subfields of $D$ which are separable over $Z(D)$.

 \item There exists a maximal subfield of $D$ which is separable over $Z(D)$.

 \item The maximal subfields of $D$ which are separable over $Z(D)$
       are exactly those subfields of $D$ which are maximal separable
       over $Z(D)$.

\end{enumerate}

{\bf\underline{Proof:}} The proof is an consequence of theorem 1
by regarding the following facts:
\begin{enumerate}

 \item Every maximal subfield and every Cartan-Subalgebra of $D$ contains the center of $D$.

 \item $D$ is central as $Z(D)$-Algebra.

 \item The Cartan-Subalgebras of $D^{\circ}$ as $K$-and $Z(D)$-Lie-Algebra are the same.$\diamond$

\end{enumerate}

An easy consequence of theorem 2 is:

\subsection*{Corollary 2}
Let $D$ be a finite-dimensional associative $K$-division algebra.

\begin{enumerate}[(i)]
 \item All Cartan-Subalgebras of $D^{\circ}$ are Lie-isomorphic and
       $ind_{Z(D)}(D)\cdot dim_K(Z(D))$-dimensional.

 \item If $Z(D)$ is separable over $K$ then the
       Cartan-Subalgebras of $D^{\circ}$ are exactly the separable
       maximal subfields of $D$.

 \item For a perfect field $K$ the Cartan-Subalgebras of
       $D^{\circ}$ are exactly the maximal subfields of $D$.$\diamond$

\end{enumerate}

We close this section by giving the following dimension-formula related to this corollary:
$ind_K(D)^2=dim_K(D)=dim_K(Z(A))\cdot dim_{Z(D)}(D)=dim_K(Z(D))\cdot ind_{Z(D)}(D)^2$.

\section*{Soluable Algebras}
In \cite{tb} and \cite{ss}T. Bauer and S. Siciliano prove for a finite-dimensional
associative unitary soluble $K$-Algebra $A$ with separable
radical factor algebra that the Cartan-Subalgebras of $A^{\circ}$
are exactly the centralizers of those subalgebras which are direct to the radical -- known as
radical complements. We will prove this theorem in a different way and analyze Cartan-Subalgebras of
soluble group algebras.\\
\\
Let $A$ be an associative $K$-algebra over a field $K$, $a$
be an algebraic element and $T$ be a subset of $A$. By $min_{a,K}$ we denote
the minimal polynomial of $a$ over $K$. Furthermore, we denote by $char(K)$ the
characteristic of $K$, by
$\langle T\rangle_K$ and $K[T]$ respectively the $K$-generating and
algebra-generating system respectively of $T$ in $A$.\\
\\
Our analysis is based on the following lemma (see for instance theorem 5.3.1
in \cite{sw1}):

\subsection*{Lemma 1}
Let $A$ be a finite-dimensional associative commutative unitary
$K$-algebra. $A$ is separable if and only if every element
of $A$ is separable over $K$.$\diamond$\\
\\
Let $A$ be an associative $K$-algebra. By $rad(A)$ we denote
the nil radical of $A$.

\subsection*{Theorem 3}
Let $A$ be a finite-dimensional associative unitary soluble
$K$-algebra with separable radical complement. The maximal
tori of $A$ are exactly the radical complements of $A$.\\
\\
\underline{\bf Proof:} '$\rightarrow :$' Let $T$ be a radical complement of $A$.
Then $T$ is a commutative separable unitary subalgebra of $A$. Using lemma 1
we conclude that $T$ is a torus of $A$. Let $S$
be a torus von $A$ which includes $T$. Again by lemma 1 the algebra $S$
is a separable $K$-subalgebra of $A$ which is direct to $rad(A)$.
Calculating the dimensions we conclude $T=S$.\\
'$\leftarrow :$' Let $T$ be a maximal torus of $A$. $T$ is
-- as a consequence of lemma 1 -- a separable subalgebra of $A$. By an enhancement
of the Wedderburn-Malcev-Conjugacy-Theorem (see
e.g. corollary 2.3.7 in \cite{sw1}) $T$ is a subalgebra of a radical complement of $A$.
The proof is completed using '$\rightarrow$'.$\diamond$\\
\\
By theorem 3, theorem 1 in \cite{ss} and the Wedderburn-Malcev-Theorem we get:

\subsection*{Theorem 4 (Bauer)}
Let $A$ be a finite-dimensional associative unitary soluble
$K$-algebra with separable radical factor algebra.

\begin{enumerate}[(i)]

 \item The Cartan-Subalgebras of $A^{\circ}$ are exactly the
       centralizers of the radical complements of $A$.

 \item The Cartan-Subalgebras of $A^{\circ}$ are conjugated
       with respect to the normal subgroup $1_A+rad(A)$ of its the group
       of units.$\diamond$

\end{enumerate}

\subsection*{Cartan-Subalgebras of soluble Group Algebras}
Let $K$ be a field and $G$ be a finite group. From 3.2.20 in
\cite{sw1} we conclude that $KG$ is soluble if and only if $G$ is
abelian or $char(K)=p$ and $G^{'}$ -- the derivation of $G$ -- is a
$p$-group. For abelian $G$ clearly $(KG)^{\circ}$
is nilpotent. Let $char(K)=p$ and $G^{'}$ be a $p$-group. As
$G^{'}$ is a normal $p$-subgroup of $G$ we conclude by Sylows-Theorems
that $G$ has exactly one (normal) $p$-Sylowsubgroup $P$.
By the Schur-Zassenhaus-Theorem there exists a complement $H$
of $P$ in $G$. Let $\alpha$ be the linearization of the canonical
group-epimorphism from $G$ onto the factor group $G/P$. Then the kernel
is given by $Kern\,\alpha = KGAug(KP)=Aug(KP)KG$ ($Aug(KP)$ is the
augmentation ideal of $KP$.). By a theorem of Wallace $Aug(KP)$ is
nilpotent and hence $Kern\,\alpha$ is nilpotent as well. The factor algebra $KG$ modulo $Kern\,\alpha$
is isomorphic to $K(G/P)$ and hence isomorphic to $KH$ as well. By Maschke's theorem $KH$
is semi-simple and hence separable using 1.9.4 in \cite{sw1}. We conclude
that $rad(KG)=KG\,Aug(KP)$ is valid and $KH$ is a
separable radical complement of $KG$.\\ By Theorem 4 all Cartan-Subalgebras of
$(KG)^{\circ}$ are conjugated by $1+rad(KG)$ to
$C_{KG}(KH)=C_{rad(KG)}(KH)\oplus KH$.\\ Using standard linear algebra we conclude that the
set $\{(a-1)h \mid 1\ne a\in P, h\in H \}$ is a $K$-basis of $rad(KG)$ which is useful
for the determination of $C_{rad(KG)}(KH)$.$\diamond$\\
\\
Let $G$ be a group, $T$ be a subset of $G$ and $a$ be
an element of finite order of $G$. By $\langle T\rangle$ we denote the subgroup of
$G$ generated by $T$ and by $o(a)$ the order of $a$ in $G$.

\subsection*{Cartan-Subalgebras of soluble Group Algebras over Dihedral Groups}
(i) Let $G$ be a group, $n\in \mathbb{N}$ and $a,b\in G$ such that
$o(a)=n$, $o(b)=2$, $G=\langle a,b\rangle$ and $a^b=a^{-1}$
are valid. If $n$ is not divisible by $2$ the derivation of $G$ is given
by $G^{'}=\langle a\rangle$. In the other case the equation
$G^{'}=\langle a^2\rangle$ is true.\\
\\
(ii) Let $K$ be a field. By 3.2.20 in
\cite{sw1} we conclude that $KG$ is soluble if and only if $G$ is
abelian or $char(K)=p$ and $G^{'}$ is a $p$-group. By (i) the following cases are to be analyzed
for the determination of Cartan-Subalgebras ($p$ is a prime not equal to $2$):\\ (a) $G$ is abelian. \\(b) $G$
is a $2$-group and $char(K)=2$.\\ (c) $n$ is a power of $p$ and $char(K)=p$.\\
(d) $\frac{n}{2}$ is a power of $p$ and $char(K)=p$.\\ Now we describe the Cartan-Subalgebras of
$(KG)^{\circ}$ for this four cases and use the conclusions of the section {\bf{Cartan-Subalgebras of soluble
Group Algebras}} for our analysis.\\
\\
(iii)(a) $G$ is abelian only for $n\in
\haken{2}$. Then $(KG)^{\circ}$ is nilpotent.\\
\\
(iii)(b) Let $G$ be a $2$-group and $char(K)=2$. Then $KG=Aug(KG)\oplus K\cdot1_G$ holds
and $(KG)^{\circ}$ is nilpotent.\\
\\
(iii)(c) Let $n$ be a power of $p$ and $char(K)=p$.
$G^{'}=\langle a\rangle$ is the $p$-Sylowsubgroup of $G$ with
complement $\langle b\rangle$. $K\langle b\rangle$ is a
radical complement of $KG$ with $K$-basis $\{1,b\}$, and the set
$\{a^s-1,(a^s-1)b \mid s\in \haken{n-1} \}$ is a $K$-basis of the
radical. The Cartan-Subalgebras are conjugated by $1+rad(KG)$
to $C_{KG}(K\langle b\rangle)=C_{rad(KG)}(K\langle
b\rangle)\oplus K\langle b\rangle$. We calculate the centralizer
of $K\langle b\rangle$ in $rad(KG)$ and show that its dimension
is $n-1$. Hence all Cartan-Subalgebras of $(KG)^{\circ}$ are $(n+1)-$ dimensional.\\ Let
$x\in rad(KG)$ like
$x=\sum\limits_{i=1}^{n-1}k_i(a^{i}-1)+\sum\limits_{j=1}^{n-1}l_j(a^j-1)b$.
We calculate: {\begin{eqnarray*} \forall y\in K\langle b\rangle: xy & =
& yx \Longleftrightarrow \\
\sum\limits_{i=1}^{n-1}k_ia^{i}b+\sum\limits_{j=1}^{n-1}l_ja^j
-\sum\limits_{i=1}^{n-1}k_iba^{i}-\sum\limits_{j=1}^{n-1}l_jba^jb
&  = & 0\Longleftrightarrow\\
\sum\limits_{i=1}^{n-1}k_i(a^{i}b-a^{-i}b)+\sum\limits_{j=1}^{n-1}l_j(a^j-a^{-j})
& = & 0\Longleftrightarrow\\
\sum\limits_{i=1}^{n-1}k_i(a^{i}b-a^{n-i}b)+\sum\limits_{j=1}^{n-1}l_j(a^j-a^{n-j})
& = & 0\Longleftrightarrow\\
\sum\limits_{i=1}^{\frac{n-1}{2}}(k_i-k_{n-i})a^{i}b+\sum\limits_{j=1}^{\frac{n-1}{2}}(l_j-l_{n-j})a^j
& = & 0\Longleftrightarrow\\ \forall i,j\in \haken{\frac{n-1}{2}}:
k_i=k_{n-i} \wedge l_j=l_{n-j} & &\\
\end{eqnarray*}}

(iii)(d) Let $\frac{n}{2}$ be a power of $p$ like $n=2\cdot p^r$ and
$char(K)=p$. $G^{'}=\langle a^2\rangle$ is the
$p$-Sylowsubgroup of $G$ with complement
$H:=\{1,b,a^{p^r},a^{p^r}b \}$. $KH$ is a radical complement in
$KG$ with $K$-basis $H$, and the set
$\{a^{2s}-1,b(a^{2s}-1),a^{p^r}a^{2s}-1,a^{p^r}ba^{2s}-1 \mid s\in
\haken{p^r-1} \}$ is a $K$-basis of the radical. The
Cartan-Subalgebras are conjugated by $1+rad(KG)$ to
$C_{KG}(KH)=C_{rad(KG)}(KH)\oplus KH$. We calculate the centralizer
of $KH$ in $rad(KG)$ and determine that its dimension is
$n-2$. Hence all Cartan-Subalgebras of
von $(KG)^{\circ}$ are of dimension $n+2$.\\ Let $x\in rad(KG)$ like
\begin{eqnarray*}
x & = &\sum\limits_{i=1}^{p^r-1}l_i(a^{2i}-1) +
\sum\limits_{i=1}^{p^r-1}m_ib(a^{2i}-1)\\ & + &
\sum\limits_{i=1}^{p^r-1}r_ia^{p^r}(a^{2i}-1) +
\sum\limits_{i=1}^{p^r-1}s_ia^{p^r}b(a^{2i}-1).
\end{eqnarray*}
$x$ centralizes $KH$ if and only if $x\circ b =0=x\circ
a^{p^r}$ is valid.\\We calculate
\begin{eqnarray*}
x\circ a^{p^r}& = &
\sum\limits_{i=1}^{p^r-1}m_i(ba^{2i}a^{p^r}-a^{p^r}ba^{2i}-ba^{p^r}+a^{p^r}b)\\
& + &
\sum\limits_{i=1}^{p^r-1}s_i(a^{p^r}ba^{2i}a^{p^r}-ba^{p^r}-a^{p^r}ba^{p^r}+b)\\
& = & 0,
\end{eqnarray*}
since $a^{p^r}$ is an involution.\\ Additionally we calculate
\begin{eqnarray*}
x\circ b & = & \sum\limits_{i=1}^{p^r-1}l_i(a^{2i}b-ba^{2i}) +
\sum\limits_{i=1}^{p^r-1}m_i(ba^{2i}b-a^{2i})\\ & + &
\sum\limits_{i=1}^{p^r-1}r_i(a^{p^r}
a^{2i}b-ba^{p^r}a^{2i}-a^{p^r}b+ba^{p^r})\\ & + &
\sum\limits_{i=1}^{p^r-1}s_i(a^{p^r}ba^{2i}b-ba^{p^r}ba^{2i}-a^{p^r}bb+ba^{p^r}b)\\
& = & \sum\limits_{i=1}^{p^r-1}l_i(a^{2i}-a^{-2i})b +
\sum\limits_{i=1}^{p^r-1}(-m_i)(a^{2i}-a^{-2i})\\ & + &
\sum\limits_{i=1}^{p^r-1}r_i(a^{2i}-a^{-2i})a^{p^r}b +
\sum\limits_{i=1}^{p^r-1}(-s_i)(a^{2i}-a^{-2i})a^{p^r}\\ & = &
\sum\limits_{i=1}^{\frac{p^r-1}{2}}(l_i-l_{p^r-i})a^{2i}b +
\sum\limits_{i=1}^{\frac{p^r-1}{2}}(m_{p^r-i}-m_i)a^{2i}\\ & + &
\sum\limits_{i=1}^{\frac{p^r-1}{2}}(r_i-r_{p^r-i})a^{2i}a^{p^r}b
+\sum\limits_{i=1}^{\frac{p^r-1}{2}}(s_{p^r-i}-s_i)a^{2i}a^{p^r}.\\
\end{eqnarray*}
Hence $x$ centralize $b$ if and only if for all
$i\in\haken{\frac{p^r-1}{2}}$ the conditions $l_i=l_{p^r-i}$,
$m_i=m_{^r-i}$, $r_i=r_{p^r-i}$ and $s_i=s_{p^r-i}$
are valid.$\diamond$

\section*{Lie-nilpotent Associative Algebras}
As we have deduced by proposition 1 the Cartan-Subalgebras of Lie-Algebras
associated to associative algebras are associative subalgebras. Hence we are interested in the
associative structure of these special subalgebras. In particular we clarify
which associative algebras are Lie-nilpotent. This topic is linked to the question
for which unitary associative algebras its group of units is nilpotent. We demonstrate are conclusions for
group algebras.\\
\\
A first result is an easy consequence of theorem 1 in \cite{ss}:

\subsection*{Corollary 1}
Let $A$ be a finite-dimensional associative unitary $K$-algebra.
$A^{\circ}$ is nilpotent if and only if every separable
element of $A$ is central.\\
\\
{\bf\underline{Proof:}} By Theorem 1 in \cite{ss} we conclude that
$A^{\circ}$ is nilpotent if and only if every maximal
torus of $A$ is central. If $a$ is a separable element of $a$ then
the subalgebra $K[a]$ is a torus (see lemma 5.2.5 in \cite{sw1}) of $A$. Hence this torus
is contained in a maximal torus of $A$.$\diamond$

\subsection*{Remark 4}
Let $A$ be a finite-dimensional associative $K$-algebra with a
central radical complement. For all $n\in \mathbb{N}$ the equation
$\underbrace{A\circ \cdots \circ A}_{n-mal} =
\underbrace{rad(A)\circ \cdots \circ rad(A)}_{n-mal}$ is valid.
As $rad(A)$ is a nilpotent associative algebra the Lie-Algebra $rad(A)^{\circ}$
is nilpotent as well. Hence $A^{\circ}$ is nilpotent with the same nilpotency class as
$rad(A)^{\circ}$.$\diamond$

\subsection*{Lemma 2}
Every finite-dimensional associative Lie-nilpotent $K$-algebra
is a soluble associative $K$-algebra.\\
\\
{\bf\underline{Proof:}} {\underline{Step 1:}} Let $A$ be a
central $K$-division algebra, $n:=ind(D)$ and $T$ be a
maximal subfield of $A$. Then $A\otimes T$ is isomorphic to $T^{n\times
n}$. By using remarks 2 and 3 we conclude $A=K\cdot 1_A$.\\
\\
{\underline{Step 2:}} Let $A$ be a $K$-division algebra. $A$ is as $Z(A)$-algebra
central and still Lie-nilpotent. By step 1 we get $A=Z(A)$.\\
\\
{\underline{Step 3:}} Let $A$ be simple. Then a
$K$-division algebra $D$ and $n\in \mathbb{N}$ exist such that $A$ is isomorphic to
$D^{n\times n}$. As $A$ is Lie-nilpotent so $D$ is Lie-nilpotent, too.
By step 2 we conclude that $D$ is a field. In addition
$A^{\circ}$ contains a subalgebra isomorphic to $gl(n,K)$. By remark
3 we get $n=1$, and hence $A$ is a field.\\
\\
{\underline{Step 4:}} Let $A$ semi-simple and therefor a direct product
of simple ideals of $A$. By step 3 every simple ideal is a field and hence $A$ is
commutative.\\
\\
{\underline{Step 5:}} By
$A^{\circ}/{rad(A)^{\circ}}=(A/rad(A))^{\circ}$ and step 4 we get that $A$ is
soluble.$\diamond$\\
\\
By lemma 2 and proposition 1 we conclude:

\subsection*{Corollary 3}
Let $A$ be a finite-dimensional associative unitary
$K$-algebra. Every Cartan-Subalgebra of $A^{\circ}$ is a
soluble associative unitary subalgebra of $A$.$\diamond$\\
\\
A detailed description of the associative structure is given by the next theorem.

\subsection*{Theorem 5}
Let $A$ be a finite-dimensional associative unitary $K$-algebra with
separable radical factor algebra. The following statements are equivalent:

\begin{enumerate}[(i)]
 \item $A^{\circ}$ is nilpotent.

 \item $A$ has a central radical complement.

 \item $A$ is soluble and has exactly one radical complement.

 \item $A$ is soluble and the set of all separable elements is a
       radical complement.

 \item $A$ is soluble and the set of all separable elements is a
       $K$-subspace.

\end{enumerate}

{\bf\underline{Proof:}} (ii) $\rightarrow$ (i): see remark 4\\
\\
(i) $\rightarrow$ (ii): Let $A^{\circ}$ be nilpotent. Using lemma 2 we get that
$A$ is soluble. Let $T$ be a radical complement of $A$. By theorem 4 we get $C_A(T)=A$. Hence
$T$ is central in $A$. By corollary 2.3.7 in \cite{sw1} we conclude (ii).\\
\\
(ii) $\rightarrow$ (iii): As the radical complement is central it is
commutative and hence $A$ is soluble. In addition
$A$ possesses exactly one radical complement .\\
\\
(iii) $\rightarrow$ (ii): Let a $A$ be soluble and it exists exactly one
radical complement of $A$. By corollary 5.1.5 in
\cite{sw1} the intersection of all radical complements is central.\\
\\
(iv) $\rightarrow$ (iii): The set of separable elements of
$A$ is invariant under conjugation by the group of units. 
Hence (iii) is a consequence of (iv).\\
\\
(iii) $\rightarrow$ (iv): Let $T$ be a radical complement of $A$
and $A$ soluble. $T$ is a commutative separable
subalgebra of $A$. By lemma 1 every
element of $T$ is separable over $K$.\\ Let $a$ be a
separable element of $A$ then $K[a]$ is due to lemma 1 a
separable subalgebra of $A$. This subalgebra is -- by using corollary 2.3.7 in
\cite{sw1} -- contained in a radical complement of $A$ and hence contained in $T$.\\
\\
(iv) $\rightarrow$ (v): This is obvious.\\
\\
(v) $\rightarrow$ (iv): Let $T$ be a radical complement of $A$.
$T$ is commutative and separable. By theorem 5.6.11 in
\cite{sw1} we conclude that every element of $T$ is separable over
$K$. The set of all separable elements is direct
to the radical. By dimension reasons we conclude that $T$ is
exactly the set of all separable elements of $A$.$\diamond$\\
\\
As a corollary we get a generalization of remark 2:

\subsection*{Corollary 2}
Let $A,B$ finite-dimensional associative $K$-algebras with
separable radical factor algebras, $A^{\circ}$ and
$B^{\circ}$ be nilpotent. Then $(A\otimes B)^{\circ}$ is nilpotent and
$A\otimes B$ has a separable radical factor algebra.\\
\\
{\bf\underline{Proof:}} Let $T$ and $S$ be radical complements of $A$ and $B$.
These complements are central by theorem 5. Theorem 2.2.9 in
\cite{sw1} shows us that $T\otimes S$ is a separable
radical complement of $A\otimes B$. As $Z(A\otimes
B)=Z(A)\otimes Z(B)$ this radical complement is central and hence
$A\otimes B$ is Lie-nilpotent using theorem 5 again.$\diamond$

\subsection*{Proposition 4}
Let $A$ be a finite-dimensional associative unitary $K$-algebra.

\begin{enumerate}[(i)]

 \item If $rad(A)$ has a central complement then $E(A)$ is nilpotent.

 \item Is $T$ a radical complement of $A$ then $1_A\in T$.

 \item Is $T$ a radical complement of $A$ and $E(T)$ central in $A$ then
       $E(A)$ is a direct product of the nilpotent normal subgroup
       $1_A+rad(A)$ and the central normal subgroup $E(T)$. In particular
       $E(A)$ is nilpotent.

\end{enumerate}

{\bf\underline{Proof:}} see 1.1.8 in \cite{sw2} and
1.10.1 in \cite{sw1}$\diamond$

\subsection*{Lemma 3}
Let $A$ be a finite-dimensional associative unitary $K$-algebra with
nilpotent group of units $E(A)$. Then $A$ is soluble.\\
\\
{\bf\underline{Proof:}} By lemma A.1.1 in \cite{sw1} we get
$E(A/rad(A))=E(A)/(1+rad(A))$ which is nilpotent as well. Let $D_1,\cdots , D_r$
be $K$-division algebras and $n_1, \cdots n_r\in \mathbb{N}$ such that
$A/rad(A)$ is isomorphic to $D_1^{n_1\times n_1} \times \cdots \times
D_r^{n_1\times n_r}$. By induction we get from 1.1.8
in \cite{sw2} that $E(A)/(1+rad(A))$ is isomorphic to $E(D_1^{n_1\times n_1})
\times \cdots \times E(D_r^{n_1\times n_r})$.
In particular we have that for all $i\in \haken{r}$ the group $E(D_i)$
is nilpotent. By a theorem of Stuth (see Corollary 5.3.1.2 in
\cite{gk}) we conclude that for all $i\in \haken{r}$ the division algebra $D_i$
is a field. In addition for all $i\in \haken{r}$
the group $GL(n_i,K)$ is nilpotent. Apart from $GL(2,2)$ and $GL(2,3)$
these groups (see page 181 in \cite{bh}) are not soluble for $n_i \ge 2$.
By page 183 in \cite{bh} the groups $GL(2,2)$ and $PSL(2,3)$ are isomorphic to
$S_3$ and $A_4$ which are not nilpotent as well. Thus we conclude $n_i=1$ for
all $i\in\haken{r}$.$\diamond$

\subsection*{Theorem 6}
Let $K$ be a field with more than 2 elements and $A$ be an
associative finite-dimensional unitary $K$-algebra with
separable radical factor algebra. The following statements are equivalent:

\begin{enumerate}[(i)]

 \item $E(A)$ is nilpotent.

 \item $A$ possesses a central radical complement.

 \item $A^{\circ}$ is nilpotent.

\end{enumerate}

{\bf\underline{Proof:}} Due to theorem 5 and proposition 4 we have to prove only
the implication (i) $\rightarrow$ (ii). Let $E(A)$ be nilpotent. Then $A$ is
soluble by lemma 3. Let $T$ be a radical complement of $A$
(Wedderburn-Malcev). By theorem 5.16 and corollary 5.18 in \cite{tb}
the group $E(C_A(T))$ is a Carter-Subgroup of $E(A)$. As a Carter-Subgroup is maximal
nilpotent we get $E(A)=E(C_A(T))$. In particular we get $1+rad(A)\subseteq
E(C_A(T)) \subseteq C_A(T)$. Hence $1_A+rad(A)$ centralizes the
radical complement $T$. By the Wedderburn-Malcev-Theorem we conclude that $A$
possesses exactly one radical complement. Due to corollary 5.1.5 in \cite{sw1} the
intersection of all radical complements of
$A$ -- which is in our case $T$ -- is central.$\diamond$

\subsection*{Theorem 7}
Let $A$ be a finite-dimensional associative unitary $K$-algebra
with separable radical factor algebra. $E(A)$ is nilpotent if and only if
for every radical complement $T$ of $A$ the group
$E(T)$ is central. In that case $E(A)$ is the direct product
of the nilpotent normal subgroup $1+rad(A)$ and the central
normal subgroup $E(T)$. Sufficient for the nilpotency of $E(A)$ is
the nilpotency of $A^{\circ}$.\\
\\
{\bf\underline{Proof:}} One implication is the content of proposition
4. Let $E(A)$ be nilpotent. Then $A$ is soluble by lemma 3.
For a radical complement $T$ of $A$ the subgroup $C_{E(A)}(E(T))$
is a Carter-Subgroup of $E(A)$ (see theorem 5.16 in \cite{tb}).
As Carter-Subgroups are maximal nilpotent we conclude $E(A)=C_{E(A)}(E(T))$.
Hence $E(T)$ is central. The proof is completed by 1.1.8
in \cite{sw2} and theorem 5.$\diamond$

\subsection*{Group Algebras}
Let $K$ be a field and $G$ be a finite group. The authors of
\cite{ab} and \cite{ik} prove that $KG$ is
Lie-nilpotent if and only if $E(KG)$ is nilpotent. This is equivalent to
$char(K)=0$ and $G$ is abelian or -- in case $char(K)=p$ -- $G^{'}$
is a $p$-subgroup of the nilpotent group $G$.\\ With respect to
theorems 6 and 7 we show that $KG$ has a separable radical factor
algebra and a central radical complement in the modular case.\\
Let $char(K)=p$ and $G^{'}$ be a $p$-subgroup of the nilpotent group $G$.
Let $P$ be the normal $p$-Sylowsubgroup of $G$ with normal complement
$N$. By $G^{'}\le P$ the normal subgroup $N$ is central in $G$. Let $\alpha$ be the
linearization of the canonical epimorphism from $G$ onto
the factor group $G/P$. Then it is well-known hat $Kern\,\alpha =
KGAug(KP)=Aug(KP)KG$ and $Aug(KP)$ is the augmentation ideal of
$KP$. By a theorem of Wallace $Aug(KP)$ is nilpotent. As a consequence
$Kern\,\alpha$ is nilpotent, too. The factor algebra
of $KG$ modulo $Kern\,\alpha$ is isomorphic to $K(G/P)$ and hence isomorphic to $KN$.
$KN$ is by theorem of Maschke semi-simple and therefor
separable by 1.9.4 in \cite{sw1}. We conclude
that $rad(KG)=KG\,Aug(KP)$ holds and that $KN$ is a separable
radical complement in $KG$. As $N$ is central in $G$ the algebra $KN$ is
a central radical complement of $KG$.$\diamond$

\section*{Simple, semi-simple and separable algebras}
We analyze Cartan-Subalgebras of Lie-Algebras
associated to simple finite-dimensional associative
$K$-algebras and reduce the analysis of semi-simple (and separable)
to their simple components. The following lemma is a version of the corresponding
one from S. Siciliano in \cite{ss}:

\subsection*{Lemma 4}
Let $A$ be a central-simple finite-dimensional associative
$K$-algebra. Is $H$ a Cartan-Subalgebra of $A^{\circ}$ then
$H$ is a $ind(A)$-dimensional self-centralizing torus
of $A$. In particular $H$ is a maximal commutative
subalgebra of $A$.\\
\\
{\bf\underline{Proof:}} By proposition 1 and lemma 2 the algebra $H$ is
a soluble unitary subalgebra of $A$. Let $n:=ind(A)$. For a
maximal subfield $T$ of $A$ the tensor product $A\otimes T$ is isomorphic
to $T^{n\times n}$. Is $F$ an algebraic closure of $T$ then $A\otimes F \cong
F^{n\times n}$ holds, and we conduct $n^2=dim_K(A)=dim_F(A\otimes F)$.
By \cite{nj} the algebra $H\otimes F$
is a Cartan-Subalgebra of the $F$-algebra $(A\otimes F)^{\circ}$ isomorphic
to $gl(n,F)$. Again by \cite{nj} all
Cartan-Subalgebras of $(F^{n\times n})^{\circ}$ are conjugate
under $GL(n,F)$ to the diagonal-matrix-algebra $D(n,F)$. Hence every
element of the $F$-algebra $H\otimes
F$ is diagonalizable and therefor every element of $H$ is separable over $K$.
In particular $H$ is semi-simple. As $H$ is soluble (lemma 2) $H$ is a torus.
$H$ is self-centralizing as $H$ is commutative
and self-normalizing: $H\subseteq C_A(H)=C_{A^{\circ}}(H)\subseteq
N_{A^{\circ}}(H)=H$. Finally for each commutative subalgebra $C$ of $A$
that contains $H$ the statement $C\subseteq C_A(C)\subseteq C_A(H)=H$ is valid.$\diamond$

\subsection*{Theorem 8}
Let $A$ be a central-simple finite-dimensional associative
$K$-algebra. The following statements are valid:

\begin{enumerate}[(i)]

 \item The maximal tori of $A$ are exactly the self-centralizing tori
       of $A$. In particular, every maximal torus of $A$ is
       a maximal commutative, separable subalgebra of $A$.

 \item The Cartan-Subalgebras of $A^{\circ}$ are exactly the maximal tori of $A$.

 \item Every Cartan-Subalgebra of $A^{\circ}$ is $ind(A)$-dimensional.

 \item All Cartan-Subalgebras of $A^{\circ}$ are isomorphic.

\end{enumerate}
{\bf\underline{Proof:}} ad(i): Let $T$ be a maximal torus of
$A$. Using theorem 1 in \cite{ss} the subalgebra $C_A(T)$ is a
Cartan-Subalgebra of $A^{\circ}$ which is by lemma 4 a
torus as well. As $T$ is commutative we conclude $T=C_A(T)$.
For a commutative subalgebra $C$ of $A$ which contains $T$ we deduct:
$C\subseteq C_A(C)\subseteq C_A(T)=T$. Finally $T$ is separable by lemma 1.\\
\\
ad(ii): This statement is a consequence of (i) and theorem 1 in \cite{ss}.\\
\\
ad(iii): see lemma 4.\\ \\ ad(iv): All Cartan-Subalgebras of
$A^{\circ}$ are abelian by (ii) and by (iii) of
the same $K$-dimension. $\diamond$\\
\\
Every simple finite-dimensional associative algebra $A$ is
central as $Z(A)$-algebra. Every Cartan-Subalgebra of
$A^{\circ}$ contains the center of $A$ and is therefor
by proposition 1 a $Z(A)$-subalgebra of $A$. Hence we conclude
by theorem 8:

\subsection*{Theorem 9}
Let $A$ be a simple finite-dimensional associative
$K$-algebra.\\ The Cartan-Subalgebras of $A^{\circ}$ are exactly
those unitary commutative subalgebras of $A$ which are maximal with
respect to that every element is separable over the
center of $A$. These subalgebras are self-centralizing and
maximal commutative.\\ In particular each Cartan-Subalgebra $T$
of $A^{\circ}$ is a direct sum of fields and
$ind_{Z(A)}(A)$-dimensional. All
Cartan-Subalgebras of $A^{\circ}$ are isomorphic.$\diamond$\\
\\
By theorems 8 and 9 we conclude:

\subsection*{Corollary 3}
Let $A$ be a simple finite-dimensional associative
$K$-algebra for which the center is separable over $K$.

\begin{enumerate}[(i)]

 \item The maximal tori are exactly the self-centralizing
       tori of $A$. In particular every maximal torus of $A$ is a
       maximal commutative, separable subalgebra von $A$.

 \item The Cartan-Subalgebras of $A^{\circ}$ are exactly
       the maximal tori of $A$.

 \item Every Cartan-Subalgebra $T$ of $A^{\circ}$ is
       $ind_{Z(A)}(A)$-dimensional.

 \item The Cartan-Subalgebras of $A^{\circ}$ are isomorphic.$\diamond$

\end{enumerate}

\subsection*{Remark 5}
We will reduce now the analysis of Cartan-Subalgebras
for semi-simple associative algebras to their simple components.\\ 
Let $A,B$ be associative $K$-algebras. For all
$a_1,a_2 \in A$ and $b_1,b_2\in B$ the rule $(a_1;b_1)\circ
(a_2;b_2)=(a_1\circ a_2;b_1\circ b_2)$ is valid. Hence we conclude
for $T\subseteq A$ and $S\subseteq B$ the equation $N_{(A\times
B)^{\circ}}(T\times S) =N_{A^{\circ}}(T)\times N_{B^{\circ}}(S)$.\\\\
(i) Let $n\in \mathbb{N}$, $A_1,\cdots , A_n$ be associative
$K$-algebras and $C_1,\cdots , C_n$ Cartan-Subalgebras of
$(A_1)^{\circ},\cdots , (A_n)^{\circ}$ then $C_1\times \cdots
\times C_n$ is a Cartan-Subalgebra of $(A_1\times \cdots
A_n)^{\circ}$.\\ \\
(ii) Let $A,B$ be associative finite-dimensional $K$-algebras and
$C$ be a Cartan-Subalgebra of $(A\times B)^{\circ}$. We define $T:=\{a \mid a\in
A, \exists b \in B: (a,b)\in C \}$ and $S:=\{b\mid b\in B, \exists
a \in A: (a,b)\in C \}$. $T$ resp. $S$ is a nilpotent
subalgebra of $A^{\circ}$ resp. $B^{\circ}$. In particular
$T\times S$ is a nilpotent subalgebra of $(A\times
B)^{\circ}$ containing (by definition) $C$ as a subalgebra. As
$C$ is maximal nilpotent we conclude $C =T\times S$. By
$T\times S=C=N_{(A\times B)^{\circ}}(C)=N_{(A\times
B)^{\circ}}(S\times T)=N_{A^{\circ}}(T)\times N_{B^{\circ}}(T)$
we deduct finally that $T$ resp. $S$ is a Cartan-Subalgebra of
$A^{\circ}$ resp. $B^{\circ}$.$\diamond$\\
\\
We deduct the following reduction-theorem which includes the description
of Cartan-Subalgebras in the semi-simple case:

\subsection*{Theorem 10}
Let $n\in \mathbb{N}$ and $A_1,\cdots , A_n$ be finite-dimensional associative
$K$-algebras. The Cartan-Subalgebras of $(A_1\times
\cdots \times A_n)^{\circ}$ are exactly the subalgebras $C_1\times \cdots
\times C_n$, whereas for every $i\in \haken{n}$ the set $C_i$ is a
Cartan-Subalgebra of $(A_i)^{\circ}$.$\diamond$\\
\\
In particular we conclude from corollary 3, theorems 10 and theorem 1 in
\cite{ss}:

\subsection*{Theorem 11}
Let $A$ be a finite-dimensional associative separable
$K$-algebra.

\begin{enumerate}[(i)]

 \item Every maximal torus of $A$ is a direct sum of maximal tori
       linked to the direct composition of $A$ into simple ideals of $A$.

 \item The maximal tori of $A$ are exactly the self-centralizing
       tori of $A$. In particular every maximal torus is a
       maximal commutative, separable subalgebra of $A$.

 \item The Cartan-Subalgebras of $A^{\circ}$ are exactly the maximal tori of $A$.

 \item All Cartan-Subalgebras of $A^{\circ}$ are isomorphic.$\diamond$

\end{enumerate}

\section*{Reduced Algebras}
For an associative $K$-algebra
$A$ we denote by $nil(A)$ the set of all nilpotent elements of $A$.
The nil radical of $A$ is always a subset of $nil(A)$.

\subsection*{Definition 2}
An associative $K$-algebra is called reduced if and only if
$rad(A)=nil(A)$ is valid.\\
\\
We will reduce the analysis of Cartan-Subalgebras of reduced associative algebras
to some special soluble subalgebras.\\
An easy observation of reduced algebras is given in the following proposition:

\subsection*{Proposition 5}
Let $A$ be a finite-dimensional associative $K$-algebra. The following statements are
equivalent:

\begin{enumerate}[(i)]

 \item $A$ is reduced.

 \item $A/rad(A)$ is reduced.

 \item $A/rad(A)$ is a direct sum of $K$-division algebras.\\
       In particular $A$ is reduced when $A$ is commutative.

 \item For every subalgebra $T$ of $A$ the condition $rad(T)=rad(A)\cap T$
        is valid.

 \item Every subalgebra of $A$ is reduced.$\diamond$

\end{enumerate}

By using lemma 1 the following remark is valid:

\subsection*{Remark 6}
Let $A$ be a finite-dimensional associative unitary $K$-algebra,
$T$ a torus of $A$ and $S:=T\oplus rad(A)$.\\ Then $T$ is a separable
complement of the radical $rad(S)=rad(A)$ of $S$ and $S$
a soluable subalgebra of $A$.$\diamond$

\subsection*{Lemma 5}
Let $A$ be a associative finite-dimensional unitary reduced
$K$-algebra.\\
The maximal elements of the set of soluble subalgebras with separable
radical factor algebra are exactly of the form $rad(A)\oplus T$, whereas
$T$ is a maximal torus of $A$.\\
\\
{\bf\underline{Proof:}} Let $T$ be a maximal torus of $A$ and
$S:=rad(A)\oplus T$. By remark 6 the set $S$ is a soluble
subalgebra of $A$ with separable radical factor algebra. Let $B$ be a
soluble subalgebra of $A$ with separable
radical factor algebra which contains $rad(A)\oplus T$. As $A$
is reduced we conclude by proposition 5 the condition $rad(B)\subseteq
rad(A)$. As a torus $T$ is by lemma 1 a separable subalgebra
of $B$ which is by corollary 2.3.7 in \cite{sw1} contained
in a radical complement $X$ of $B$. $X$ is using lemma 1 a torus
of $A$ as well, and we deduct $T=X$ from the maximality of $T$. Therefor
we proved $B=rad(A)\oplus T$.\\ For the other implication let $B$
be a maximal element of the set of soluble subalgebras with separable
radical factor algebra and $T$ be a radical complement
of $B$. The subalgebra $T$ is by lemma 1 a torus of $A$, and by
proposition 5 we conclude $rad(B)\subseteq rad(A)$. We define
$S:=rad(A)\oplus T$. Using remark 6 and the maximality of $B$
we conclude $B=S$ and $rad(B)=rad(A)$. Suppose $R$ is a $T$ containing
torus not equal to $T$ then the subalgebra $rad(A)\oplus R$ is
containing $B$ and not equal to $B$. With remark 6 we deduct a
contradiction to the maximality of $B$.$\diamond$\\
\\
By lemma 5 we get the following corollary:

\subsection*{Corollary 4}
Let $A$ be a associative finite-dimensional unitary reduced
$K$-algebra over a perfect field $K$. The maximal soluble
subalgebras of $A$ are exactly of the form $T\oplus rad(A)$ whereas
$T$ is a maximal torus of $A$.$\diamond$

\subsection*{Remark 7}
Let $A$ be a associative $K$-algebra, $I$ an ideal and $T$ a
subalgebra of $A$ such that $A$ is the internal direct sum of $I$
and $T$. For every subset $X$ of $T$ the equation
$C_A(X)=C_I(X)\oplus C_T(X)$ is valid.$\diamond$

\subsection*{Theorem 12}
Let $A$ be a finite-dimensional associative unitary reduced
$K$-algebra with separable radical factor algebra and
$\mathcal{S}(A)$ the set of soluble subalgebras with
separable radical factor algebra of $A$.\\ The
Cartan-Subalgebras of $A^{\circ}$ are exactly
the Cartan-Subalgebras of Lie-Algebras associated to maximal
elements of $\mathcal{S}(A)$.\\
\\
{\bf\underline{Proof:}} Let $H$ be a Cartan-Subalgebra of
$A^{\circ}$. By theorem 1 in \cite{ss} there exists a
maximal torus $T$ of $A$ such that $H=C_A(T)$ is valid. We
define $S:=rad(A)\oplus T$. By lemma 5 the set $S$ is a maximal
element of $\mathcal{S}(A)$. Using 2.3.7 in \cite{sw1} and lemma 1
there exists a radical complement $C$ of $A$ containing $T$. Thus
$T$ is a maximal torus of $C$ which is by theorem 11
in $C$ self-centralizing. By remark 6 we conclude
$H=C_A(T)=C_{rad(A)}(T)\oplus T$. $T$ is of course
a maximal torus of $S$. By theorem 1 in \cite{ss} the subalgebra $C_S(T)$
is a Cartan-Subalgebra of $S^{\circ}$. As $T$ is commutative we get
by remark 7 the statement $C_S(T)=C_{rad(A)}(T)\oplus
T=C_A(T)=H$. Hence $H$ is a Cartan-Subalgebra von $S$.\\ Conversely
let $S$ be a maximal element of $\mathcal{S}(A)$. By lemma 5
there exists a maximal torus $T$ of $A$ such that
$S=rad(A)\oplus T$ is valid. Every Cartan-Subalgebra of $S$ is
due to theorem 4 the centralizer of a radical complement of $S$.
These are by the Wedderburn-Malcev-Theorem conjugates of $T$ with respect
to $1+rad(S)=1+rad(A)$. Now let $H$ be a Cartan-Subalgebra of $S$. We can assume $H=C_S(T)$.
By remark 6 and the commutativity of $T$ we deduct
$H=C_S(T)=C_{rad(A)}(T)\oplus T$. Using theorem 1 in \cite{ss} the subalgebra
$C_A(T)$ is a Cartan-Subalgebra of $A^{\circ}$. Let $C$ be a
radical complement of $A$ containing $T$ (see 2.3.7 in
\cite{sw1} and lemma 1). By remark 6 and theorem 11 we conclude
$C_A(T)=C_{rad(A)}(T)\oplus C_{C}(T)=C_{rad(A)}(T)\oplus
T=C_S(T)=H$.$\diamond$\\
\\
By theorem 12 we conclude the following corollary:

\subsection*{Corollary 5}
Let $A$ be finite-dimensional associative unitary reduced
$K$-algebra over a perfect field $K$. The Cartan-Subalgebras of $A^{\circ}$ are exactly the Cartan-Subalgebras
of Lie-Algebras associated to maximal soluble subalgebras of $A$.$\diamond$

\section*{Associative Algebras}

\subsection*{Theorem 13}
Let $A$ be a associative finite-dimensional unitary $K$-algebra
with separable radical factor algebra and $H$ be a subset of
$A$. The following statements are equivalent:

\begin{enumerate}[(i)]

 \item $H$ is a Cartan-Subalgebra of $A^{\circ}$.

 \item It exists a radical complement $C$ of $A$ and
       a Cartan-Subalgebra $T$ of $C^{\circ}$ such that
       $H=C_A(T)$.

\end{enumerate}

{\bf\underline{Proof:}} (i) $\rightarrow$ (ii): Let $H$ be a
Cartan-Subalgebra of $A^{\circ}$. By theorem 1 in \cite{ss}
there exists a maximal torus $T$ of $A$ such that $H=C_A(T)$
is valid. By 2.3.7 in \cite{sw1} and lemma 1 the torus $T$ is contained
in a radical complement $C$ of $A$. $T$ is a maximal
torus of $C$ as well. By theorem 11 the separable subalgebra $T$
is a Cartan-Subalgebra of $C^{\circ}$.\\
\\
(ii) $\rightarrow$ (i): Let $C$ be a radical complement of $A$ and
$T$ be a Cartan-Subalgebra of $C^{\circ}$. By theorem 11 the subalgebra $T$
is a maximal torus of $C$ such that $T=C_C(T)$. By remark
6 we conclude $C_A(T)=C_C(T)\oplus C_{rad(A)}(T)=T\oplus
C_{rad(A)}(T)$. Obviously $T$ is a central radical complement
in $C_A(T)$. Using remark 4 we deduct that $C_A(T)$ is Lie-nilpotent. $T$
is contained in a maximal torus $S$ of $A$. By theorem 1 in
\cite{ss} the subalgebra $C_A(S)$ is a Cartan-Subalgebra of $A^{\circ}$.
$C_A(S)$ is maximal nilpotent and contained in $C_A(T)$. Hence
we get $C_A(T)=C_A(S)$ and $C_A(T)$ is a
Cartan-Subalgebra of $A^{\circ}$.$\diamond$\\
\\
By theorems 11 and 13 we conclude the following {\bf strategy} for calculating
a Cartan-Subalgebra of Lie-Algebras associated to associative
unitary finite-dimensional algebras with separable
radical factor algebra:\\ \\(1) Determine a maximal tori of the
radical complements. These are exactly the self-centralizing tori of the
radical complements. To determine a self-centralizing torus in a
radical complement $C$ begin with a (as big as possible) torus $T$.
Calculate the centralizer of the torus in $C$ and search for a separable
element $t$ of this centralizer not contained in $T$. The subalgebra $K[T,t]$
is a torus as well which contains $T$. Repeat this approach until no such
separable element can be determined. \\ \\(2) For determining a Cartan-Subalgebras of
$A^{\circ}$ calculate for a maximal torus of (1) the centralizer in $A$ (theorem 13).
This centralizer is by remark 6 and theorem 11 exactly $C_{rad(A)}(T)\oplus T$.\\ \\
This approach is demonstrated now by giving an alternative proof of theorem 4 and
by determining a Cartan-Subalgebra of $(KD_{2n})^{\circ}$.

\subsection*{Corollary 5}
Let $A$ be a finite-dimensional associative unitary soluble
$K$-algebra with separable radical factor algebra. The
Cartan-Subalgebras of $A^{\circ}$ are exactly the centralizers
of the radical complements of $A$.\\
\\
{\bf\underline{Proof:}} Let $C$ be a radical complement of $A$.
By using lemma 1 the subset $C$ is a torus. The only maximal torus of $C$
is $C$ itself.$\diamond$

\subsection*{Group algebras of Dihedral Groups}
Let $n\in \mathbb{N}$, $G:=D_{2n}$ a Dihedral Group and $K$ be a
field with $char(K)=p$. There exists $a,b\in G$ such that $o(a)=n$,
$o(b)=2$, $G=\langle a,b\rangle$, $a^b=a^{-1}$ and
$G=\{1,a,...,a^{n-1},b,ab,...,a^{n-1}b\}$ are valid. Is $p>0$ not a factor
of $n$ or $p=0$ then -- by a theorem of Maschke -- the subalgebra $K\langle
a\rangle$ is semi-simple (and commutative). By 1.9.4.2 in \cite{sw1}
the algebra $K\langle a\rangle$ is separable over $K$, and by lemma 1
we conclude that $K\langle a\rangle$ is a torus of $KG$.
We analyze the centralizer of $K\langle a\rangle$ in $KG$.
We have $KG=K\langle a\rangle\oplus \langle
b,ab,...,a^{n-1}b\rangle_K$, and $K\langle a\rangle$ centralizes
$K\langle a\rangle$. Let $x\in \langle
b,ab,...,a^{n-1}b\rangle_K$, like $x=\sum\limits_{i=0}^{n-1}
k_ia^{i}b$. $x$ centralizes $K\langle a\rangle$ if and only if
$\sum\limits_{i=0}^{n-1} k_ia^{i}b^{a}=\sum\limits_{i=0}^{n-1}
k_ia^{i}b$ is valid. Because of $b^{a}=a^{-2}b$ this is equivalent
to $a^{-2}\, (\sum\limits_{i=0}^{n-1}
k_ia^{i})=\sum\limits_{i=0}^{n-1} k_ia^{i}$. If $2$ is not a factor
of $n$ (hence $\langle a^{-2}\rangle=\langle a\rangle$ is valid) then
$x$ centralizes $K\langle a\rangle$ if and only if $K\langle
a\rangle$ acts trivial on $\sum\limits_{i=0}^{n-1} k_ia^{i}\in
K\langle a\rangle$. This is equivalent to $\sum\limits_{i=0}^{n-1}
k_ia^{i}\in \langle\sum\limits_{i=0}^{n-1} a^{i}\rangle_K$.
We conclude $C_{KG}(K\langle a\rangle)=K\langle
a\rangle\oplus \langle b+ab+...+a^{n-1}b\rangle_K=K\langle
a\rangle\oplus \langle\sum\limits_{g\in G} g\rangle_K$.\\ \\ Is
$p$ a factor of the order of $G$ -- hence we have
$p=2$ under our conditions -- then the one-dimensional $K$-subspace
$\langle\sum\limits_{g\in G} g\rangle_K$ is a zero-ideal contained in
nil radical of $KG$. We assume that $KG/rad(KG)$ is separable.
The torus $K\langle a\rangle$ is contained by lemma 1 and 2.3.7 in
\cite{sw1} in a radical complement $C$ of $KG$. By remark 6 we have
$C_{KG}(K\langle a\rangle)=C_{rad(KG)}(K\langle a\rangle)\oplus
C_C(K\langle a\rangle)$. In addition the statements $C_{KG}(K\langle
a\rangle)=K\langle a\rangle\oplus \langle\sum\limits_{g\in G}
g\rangle_K$, $K\langle a\rangle\subseteq C_C(K\langle
a\rangle)$ and $\langle\sum\limits_{g\in G} g\rangle_K\subseteq
C_{rad(KG}(K\langle a\rangle)$ are valid. Thus $K\langle
a\rangle$ self-centralizing in $C$. By the theorems 11 and 13 the set
$C_{KG}(K\langle a\rangle)=K\langle a\rangle\oplus
\langle\sum\limits_{g\in G} g\rangle_K$ is a $(n+1)$-dimensional
Cartan-Subalgebra of $(KG)^{\circ}$.\\ \\ Let $p$ not be a factor
of the order of the group $G$. The element $s:=\sum\limits_{g\in G}
g$ is diagonalizable (and hence separable over $K$)
because of $s^2=\, \mid G \mid \, s$. Using our strategy
$C_{KG}(K\langle a\rangle)$ is a torus of $A$ which is
obviously self-centralizing in $KG$. Hence $C_{KG}(K\langle a\rangle)$
is again a Cartan-Subalgebra of $(KG)^{\circ}$.$\diamond$

\end{document}